
%
%
%
%
%
%
\magnification=\magstephalf      
%
%
\vsize=7.5truein                 
\hsize=5.2truein                 
\newskip\stdskip                 
\stdskip=6pt plus3pt minus3pt    
\medskipamount=\stdskip          
\parindent=0pt                   
\parskip=\stdskip                
\abovedisplayskip=\stdskip       
\belowdisplayskip=\stdskip       
\mathsurround=0.75pt             
\overfullrule=0pt                
%
%
\def\ppar{\par\goodbreak\vskip 8pt plus 4pt minus 4pt}     
%
%
\def\stdspace{\hskip 0.75em plus 0.15em\ignorespaces}
\let\qua\stdspace 
%
%
%
%
%
%
%
\def\hexnumber#1{\ifcase#1 0\or 1\or 2\or 3\or 4\or 5\or 6\or 7\or 8\or
 9\or A\or B\or C\or D\or E\or F\fi}
%
%
\font\thirtnmsa=msam10 scaled 1315    
\font\tenmsa=msam10          \font\ninemsa=msam9
\font\sevenmsa=msam7         \font\sixmsa=msam6
\font\fivemsa=msam5
%
%
\newfam\msafam                  \textfont\msafam=\tenmsa
\scriptfont\msafam=\sevenmsa    \scriptscriptfont\msafam=\fivemsa
\edef\hexa{\hexnumber\msafam}        
\def\msa{\fam\msafam\tenmsa}         
%
%
\font\thirtnmsb=msbm10 scaled 1315   
\font\tenmsb=msbm10      \font\ninemsb=msbm9
\font\sevenmsb=msbm7     \font\sixmsb=msbm6
\font\fivemsb=msbm5
%
\newfam\msbfam                   \textfont\msbfam=\tenmsb       
\scriptfont\msbfam=\sevenmsb     \scriptscriptfont\msbfam=\fivemsb
\edef\hexb{\hexnumber\msbfam}    
\def\msb{\fam\msbfam\tenmsb}     
%
%
\font\thirtneufm=eufm10 scaled 1315   
\font\teneufm=eufm10                 \font\nineeufm=eufm9
\font\seveneufm=eufm7                \font\sixeufm=eufm6
\font\fiveeufm=eufm5
%
\newfam\eufmfam                    \textfont\eufmfam=\teneufm
\scriptfont\eufmfam=\seveneufm     \scriptscriptfont\eufmfam=\fiveeufm
\edef\hexf{\hexnumber\eufmfam}      
\def\frak{\fam\eufmfam\teneufm}     
%
%
%
\font\thirtnrm=cmr10 scaled 1315    
\font\ninerm=cmr9                   \font\sixrm=cmr6   
%
\font\thirtni=cmmi10 scaled 1315    
\font\ninei=cmmi9                   \font\sixi=cmmi6  
%
\font\thirtnsy=cmsy10 scaled 1315   
\font\ninesy=cmsy9                  \font\sixsy=cmsy6  
%
\font\thirtnbf=cmbx10 scaled 1315   
\font\ninebf=cmbx9                  \font\sixbf=cmbx6  
%
%
\font\thirtnex=cmex10 scaled 1315   
\font\nineex=cmex9                  
%
%
\font\thirtnit=cmti10 scaled 1315  
\font\nineit=cmti9                  
%
\font\thirtnsl=cmsl10 scaled 1315  
\font\ninesl=cmsl9                  
%
\font\thirtntt=cmtt10 scaled 1315  
\font\ninett=cmtt9                  
%
%
%
%
\def\small{%
%
%
\textfont0=\ninerm \scriptfont0=\sixrm \scriptscriptfont0=\fiverm
\def\rm{\fam0\ninerm}
%
%
\textfont1=\ninei \scriptfont1=\sixi \scriptscriptfont1=\fivei
%
%
\textfont2=\ninesy \scriptfont2=\sixsy \scriptscriptfont2=\fivesy
%
%
\textfont3=\nineex \scriptfont3=\nineex \scriptscriptfont3=\nineex
%
%
\textfont\bffam=\ninebf \scriptfont\bffam=\sixbf
\scriptscriptfont\bffam=\fivebf \def\bf{\fam\bffam\ninebf}%
%
%
\textfont\itfam=\nineit \def\it{\fam\itfam\nineit}%
\textfont\slfam=\ninesl \def\sl{\fam\slfam\ninesl}%
\textfont\ttfam=\ninett \def\tt{\fam\ttfam\ninett}%
%
%
%
\textfont\msafam=\ninemsa \scriptfont\msafam=\sixmsa
\scriptscriptfont\msafam=\fivemsa \def\msa{\fam\msafam\ninemsa}%
%
%
\textfont\msbfam=\ninemsb \scriptfont\msbfam=\sixmsb
\scriptscriptfont\msbfam=\fivemsb \def\msb{\fam\msbfam\ninemsb}%
%
%
\textfont\eufmfam=\nineeufm  \scriptfont\eufmfam=\sixeufm
\scriptscriptfont\eufmfam=\fiveeufm \def\frak{\fam\eufmfam\nineeufm}%
%
%
%
\normalbaselineskip=11pt%
\setbox\strutbox=\hbox{\vrule height8pt depth3pt width0pt}%
%
%
\normalbaselines\rm
%
%
\stdskip=4pt plus2pt minus2pt    
\medskipamount=\stdskip          
\parskip=\stdskip                
\abovedisplayskip=\stdskip       
\belowdisplayskip=\stdskip       
\def\ppar{\par\goodbreak\vskip 6pt plus 3pt minus 3pt}%
%
%
\def\section##1{\global\advance\sectionnumber by 1
\vskip-\lastskip\penalty-800\vskip 20pt plus10pt minus5pt 
\egroup{\bf\number\sectionnumber\quad##1}\bgroup\small         
\vskip 6pt plus3pt minus3pt
\nobreak\resultnumber=1}
}    
%
\def\beginsmall{\bgroup\small}
\let\endsmall\egroup
%
%
%
%
\def\large{%
\textfont0=\thirtnrm \scriptfont0=\ninerm \scriptscriptfont0=\sevenrm
\def\rm{\fam0\thirtnrm}%
\textfont1=\thirtni \scriptfont1=\ninei \scriptscriptfont1=\seveni
\textfont2=\thirtnsy \scriptfont2=\ninesy \scriptscriptfont2=\sevensy
\textfont3=\thirtnex \scriptfont3=\thirtnex \scriptscriptfont3=\thirtnex
\textfont\bffam=\thirtnbf \scriptfont\bffam=\ninebf
\scriptscriptfont\bffam=\sevenbf \def\bf{\fam\bffam\thirtnbf}%
\textfont\itfam=\thirtnit \def\it{\fam\itfam\thirtnit}%
\textfont\slfam=\thirtnsl \def\sl{\fam\slfam\thirtnsl}%
\textfont\ttfam=\thirtntt \def\tt{\fam\ttfam\thirtntt}%
\textfont\msafam=\thirtnmsa \scriptfont\msafam=\ninemsa
\scriptscriptfont\msafam=\sevenmsa \def\msa{\fam\msafam\thirtnmsa}%
\textfont\msbfam=\thirtnmsb \scriptfont\msbfam=\ninemsb
\scriptscriptfont\msbfam=\sevenmsb \def\msb{\fam\msbfam\thirtnmsb}%
\textfont\eufmfam=\thirtneufm  \scriptfont\eufmfam=\nineeufm
\scriptscriptfont\eufmfam=\seveneufm \def\frak{\fam\eufmfam\teneufm}%
\normalbaselineskip=16pt%
\setbox\strutbox=\hbox{\vrule height11.5pt depth4.5pt width0pt}%
\normalbaselines\rm}%
\let\Large\large   
%
\def\Bbb#1{{\msb#1}}

%

%
\mathchardef\plussquare="0\hexa01
\mathchardef\nge="3\hexb0B
\mathchardef\maltesecross="0\hexa7A
\mathchardef\del="0\hexf01
%
%
%
%
\font\sc=cmcsc10
%
%
%
%
\def\sqr#1#2{{\vcenter{\vbox{\hrule  height.#2truept
	\hbox{\vrule width.#2truept height#1truept 
	\kern#1truept \vrule width.#2truept}
	\hrule height.#2truept}}}}
\def\sq{\sqr55}    
%
%
%
%
\newcount\sectionnumber            
\newcount\resultnumber             
\sectionnumber=0\resultnumber=1    
%
%
%
\def\section#1{\global\advance\sectionnumber by 1
\xdef\nextkey{\number\sectionnumber}
\vskip-\lastskip\penalty-800\vskip 20pt plus10pt minus5pt 
{\large\bf\number\sectionnumber\quad#1}         
\vskip 8pt plus4pt minus4pt
\nobreak\resultnumber=1}      
%
%
%
%
%
         
%
%
%
%

%
\def\proc#1{\xdef\nextkey{\number\sectionnumber.\number\resultnumber}%
\vskip-\lastskip\ppar\bf%
\noindent#1\ \number\sectionnumber.\number\resultnumber
\stdspace\sl\global\advance\resultnumber by 1\ignorespaces}
\def\endproc{\rm\ppar} 
%
%
\def\prf{\vskip-\lastskip\ppar\noindent{\bf Proof}%
\stdspace\rm}                            
\def\endprf{\unskip\stdspace\hbox{}
\hfill$\sq$\par\medskip}                 
%
%
%
%
%
%
%
%
\def\proclaim#1{\vskip-\lastskip\ppar\bf%
\noindent#1\stdspace\sl\ignorespaces} 

%
%
%
%

%
%
%
%
%
%
\def\label{\xdef\nextkey{\number\sectionnumber.\number\resultnumber}%
\number\sectionnumber.\number\resultnumber
\global\advance\resultnumber by 1}
%
%
%
%
%
%
%
%
%
%
%
%
%
%
%
%
\newcount\refnumber              
\refnumber=1                     
\long\def\reflist#1\endreflist{%
\long\def\thereflist{#1}{\def\refkey##1##2\par{\xdef##1{\number\refnumber}%
\global\advance\refnumber by 1}%
\def\key##1##2\par{\expandafter\xdef%
\csname##1\endcsname{\number\refnumber}%
\global\advance\refnumber by 1}#1\par}}
\long\def\references{%
\penalty-800\vskip-\lastskip\vskip 15pt plus10pt minus5pt 
{\large\bf References}\ppar 
{\leftskip=25pt\frenchspacing    
\small\parskip=3pt plus2pt       
\def\refkey##1##2\par{\noindent  
\llap{[##1]\stdspace}\ignorespaces##2\par}         
\def\key##1##2\par{\noindent  
\llap{[\ref{##1}]\stdspace}\ignorespaces##2\par}  
\def\,{\thinspace}\thereflist\par}}
%
%
%
\newcount\footnotenumber         
\footnotenumber=1                
\def\fnote#1{\xdef\nextkey{\number\footnotenumber}%
{\small\ifnum\footnotenumber>9\parindent=14pt%
\else\parindent=10pt\fi\footnote{$^{\number\footnotenumber}$}%
{\hglue-5pt#1}\global\advance\footnotenumber by 1}}
%
%
%
%
%
%
%
\newcount\figurenumber          
\figurenumber=1                 
\def\caption#1{\xdef\nextkey{\number\figurenumber}%
\cl{\small Figure \number\figurenumber: #1}%
\global\advance\figurenumber by 1}
\def\figurelabel{\xdef\nextkey{\number\figurenumber}%
\cl{\small Figure \number\figurenumber}%
\global\advance\figurenumber by 1}
\long\def\figure#1\endfigure{{\xdef\nextkey{\number\figurenumber}%
\let\captiontext\relax\def\caption##1{\xdef\captiontext{##1}}%
\midinsert\cl{\ignorespaces#1\unskip\unskip\unskip\unskip}\vglue6pt\cl{\small 
Figure \number\figurenumber\ifx\captiontext\relax\else: \captiontext
\fi}\endinsert\global\advance\figurenumber by 1}}
%
%
%
%
%
%
%
\def\nextkey{??}   
%
\def\key#1{\expandafter\xdef\csname #1\endcsname{\nextkey}}
\def\ref#1{\expandafter\ifx\csname #1\endcsname\relax
\immediate\write16{Reference {#1} undefined}??\else
\csname #1\endcsname\fi}
%
%
%
%
%
%
%
\newread\gtinfile
\newwrite\gtreffile
\def\useforwardrefs{
\openin\gtinfile\jobname.ref
\ifeof\gtinfile
\closein\gtinfile
\immediate\write16{No file \jobname.ref}
\else
\closein\gtinfile
\input \jobname.ref
\fi
\immediate\openout\gtreffile \jobname.ref
%
%
\def\key##1{{\def\\{\noexpand}%
\expandafter\xdef\csname ##1\endcsname{\nextkey}%
\immediate\write\gtreffile{\\\expandafter\\\def\\\csname ##1\\\endcsname%
{\nextkey}}}}
%
%
\long\def\reflist##1\endreflist{%
\long\def\thereflist{##1}{\def\refkey####1####2\par{\xdef####1{%
\number\refnumber}{\def\\{\noexpand}\immediate\write\gtreffile
{\\\def\\####1{\number\refnumber}}}\global\advance\refnumber by 1}%
\def\key####1####2\par{\expandafter\xdef%
\csname####1\endcsname{\number\refnumber}%
{\def\\{\noexpand}\immediate\write\gtreffile
{\\\expandafter\\\def\\\csname ####1\\\endcsname{\number\refnumber}}}
\global\advance\refnumber by 1}##1\par}}
\long\def\biblio##1\endbiblio{\reflist##1\endreflist\references}%
%
%
\def\numkey##1{{\def\\{\noexpand}%
\xdef##1{\number\sectionnumber.\number\resultnumber}
\immediate\write\gtreffile{\\\def\\##1%
{\number\sectionnumber.\number\resultnumber}}}}
\def\seckey##1{{\def\\{\noexpand}\xdef##1{\number\sectionnumber}
\immediate\write\gtreffile{\\\def\\##1{\number\sectionnumber}}}}
\def\figkey##1{\xdef##1{\number\figurenumber}%
{\def\\{\noexpand}\immediate\write\gtreffile%
{\\\def\\##1{\number\figurenumber}}}
\number\figurenumber\global\advance\figurenumber by 1}
}   
%
%
%
%
\def\figkey#1{\xdef#1{\number\figurenumber}%
\number\figurenumber\global\advance\figurenumber by 1}
\def\fig#1#2\endfig{%
\midinsert\cl{#2}\vglue6pt\cl{\small Figure #1}\endinsert}
\def\newfig{\number\figurenumber\global\advance\figurenumber by 1}
\def\numkey#1{\xdef#1{\number\sectionnumber.\number\resultnumber}}
\def\seckey#1{\xdef#1{\number\sectionnumber}}
%
%
%
%
%
%
%
%
%
\def\verb{\catcode`\"=\active}       
\def\brev{\catcode`\"=12}            
\brev                                
\verb                                
{\obeyspaces\gdef {\ }}              
{\catcode`\`=\active\gdef`{\relax\lq}}
\def"{%
\begingroup\baselineskip=12pt\def\par{\leavevmode\endgraf}%
\tt\obeylines\obeyspaces\parskip=0pt\parindent=0pt%
\catcode`\$=12\catcode`\&=12\catcode`\^=12\catcode`\#=12%
\catcode`\_=12\catcode`\~=12%
\catcode`\{=12\catcode`\}=12\catcode`\%=12\catcode`\\=12%
\catcode`\`=\active\let"\endgroup}
\brev      
%
%
%
%
%
%
\def\items{\par\leftskip = 25pt}           
\def\enditems{\par\leftskip = 0pt}         
\def\item#1{\par\leavevmode\llap{#1\stdspace}%
\ignorespaces}                             
%
%

%
%
\def\co{\colon\thinspace}    
\def\np{\vfil\eject}         
\def\nl{\hfil\break}         
\def\cl{\centerline}         
\def\gt{{\mathsurround=0pt\it $\cal G\mskip-2mu$eometry \&\ 
$\cal T\!\!$opology}}        
\def\agt{{\mathsurround=0pt\it$\cal A\mskip-.7mu$lgebraic \&\ 
$\cal G\mskip-2mu$eometric $\cal T\!\!$opology}}  
%
%
%

%
%
%
%
%
\def\title#1{\def\thetitle{#1}}

\def\author#1{\edef\previousauthors{\theauthors}
 \ifx\theauthors\relax\def\theauthors{#1}\else
 \def\theauthors{\previousauthors\par#1}\fi}

%
\def\address#1{\edef\previousaddresses{\theaddress}
 \ifx\theaddress\relax\def\theaddress{#1}\else
 \def\theaddress{\previousaddresses\par\vskip 2pt\par#1}\fi}
\def\secondaddress#1{\edef\previousaddresses{\theaddress}
 \ifx\theaddress\relax\def\theaddress{#1}\else
 \def\theaddress{\previousaddresses\par{\rm and}\par#1}\fi}   

\def\email#1{\edef\previousemails{\theemail}
 \ifx\theemail\relax\def\theemail{#1}\else
 \def\theemail{\previousemails\hskip 0.75em\relax#1}\fi}
\def\secondemail#1{\edef\previousemails{\theemail}
 \ifx\theemail\relax\def\theemail{#1}\else
 \def\theemail{\previousemails\hskip 0.75em{\rm and}\hskip 0.75em
 \relax#1}\fi}
\def\url#1{\edef\previousurls{\theurl}
 \ifx\theurl\relax\def\theurl{#1}\else
 \def\theurl{\previousurls\hskip 0.75em\relax#1}\fi}
\def\secondurl#1{\edef\previousurls{\theurl}
 \ifx\theurl\relax\def\theurl{#1}\else
 \def\theurl{\previousurls\hskip 0.75em{\rm and}\hskip 0.75em
 \relax#1}\fi}
\long\def\abstract#1\endabstract{\long\def\theabstract{#1}}
\def\primaryclass#1{\def\theprimaryclass{#1}}
\def\secondaryclass#1{\def\thesecondaryclass{#1}}
\def\keywords#1{\def\thekeywords{#1}}
%
%
\let\\\par\let\thetitle\relax\let\theshorttitle\relax
\let\theauthors\relax\let\theshortauthors\relax
\let\theaddress\relax\let\theshortaddress\relax
\let\theemail\relax\let\theurl\relax
\let\theabstract\relax\let\theprimaryclass\relax
\let\thesecondaryclass\relax\let\thekeywords\relax
%
%
%
%
\long\def\maketitlepage{    

\vglue 0.2truein   

%
{\parskip=0pt\leftskip 0pt plus 1fil\def\\{\par\smallskip}{\large
\bf\thetitle}\par\medskip}   

\vglue 0.15truein 

%
{\parskip=0pt\leftskip 0pt plus 1fil\def\\{\par}{\sc\theauthors}
\par\medskip}%
 
\vglue 0.1truein 

%
{\small\parskip=0pt
{\leftskip 0pt plus 1fil\def\\{\par}{\sl\theaddress}\par}
\ifx\theemail\relax\else  
\vglue 5pt \def\\{\stdspace{\rm and}\stdspace} 
\cl{Email:\stdspace\tt\theemail}\fi
\ifx\theurl\relax\else    
\vglue 5pt \def\\{\stdspace{\rm and}\stdspace} 
\cl{URL:\stdspace\tt\theurl}\fi\par}

\vglue 7pt 

{\bf Abstract}

\vglue 5pt

\theabstract

\vglue 7pt 

{\bf AMS Classification numbers}\quad Primary:\quad \theprimaryclass\par

Secondary:\quad \thesecondaryclass

\vglue 5pt 

{\bf Keywords:}\quad \thekeywords

\np  

}    
%
%
\long\def\makeshorttitle{    


%
{\parskip=0pt\leftskip 0pt plus 1fil\def\\{\par\smallskip}{\large
\bf\thetitle}\par\medskip}   

\vglue 0.05truein 

%
{\parskip=0pt\leftskip 0pt plus 1fil\def\\{\par}{\sc\theauthors}
\par\medskip}%
 
\vglue 0.03truein 

%
{\small\parskip=0pt
{\leftskip 0pt plus 1fil\def\\{\par}{\sl\ifx\theshortaddress\relax
\theaddress\else\theshortaddress\fi}\par}
\ifx\theemail\relax\else  
\vglue 5pt \def\\{\stdspace{\rm and}\stdspace} 
\cl{Email:\stdspace\tt\theemail}\fi
\ifx\theurl\relax\else    
\vglue 5pt \def\\{\stdspace{\rm and}\stdspace} 
\cl{URL:\stdspace\tt\theurl}\fi\par}

\vglue 10pt 


{\small\leftskip 25pt\rightskip 25pt{\bf Abstract}\stdspace\theabstract

{\bf AMS Classification}\stdspace\theprimaryclass
\ifx\thesecondaryclass\relax\else; \thesecondaryclass\fi\par
{\bf Keywords}\stdspace \thekeywords\par}
\vglue 7pt
}    
\let\maketitle\makeshorttitle        
%
%

\def\volumenumber#1{\def\thevolumenumber{#1}}
\def\volumeyear#1{\def\thevolumeyear{#1}}
\def\pagenumbers#1#2{\def\startpage{#1}\def\finishpage{#2}}
\def\published#1{\def\publishdate{#1}}
\def\received#1{\def\receiveddate{#1}}
\def\revised#1{\def\reviseddate{#1}}
\let\reviseddate\relax
\volumenumber{X}
\volumeyear{20XX}
\pagenumbers{1}{XXX}
\published{XX Xxxember 20XX}

\long\def\makeagttitle{   
\agt\hfill      
\hbox to 60truept{\vbox to 0pt{\vglue -14truept{\bf [Logo here]}\vss}\hss}
\break
{\small Volume \thevolumenumber\ (\thevolumeyear)
\startpage--\finishpage\nl
Published: \publishdate}

\vglue .2truein

{\parskip=0pt\leftskip 0pt plus 1fil\def\\{\par\smallskip}{\large
\bf\thetitle}\par\medskip}   
\vglue 0.05truein 

%
{\parskip=0pt\leftskip 0pt plus 1fil\def\\{\par}{\sc\theauthors}
\par\medskip}%
 
\vglue 0.03truein 


{\small\leftskip 25truept\rightskip 25truept{\bf Abstract}\stdspace\theabstract

{\bf AMS Classification}\stdspace\theprimaryclass
\ifx\thesecondaryclass\relax\else; \thesecondaryclass\fi\par
{\bf Keywords}\stdspace \thekeywords\par}\vglue 7truept

}   


\def\Addresses{\bigskip
{\small \parskip 0pt \leftskip 0pt \rightskip 0pt plus 1fil \def\\{\par}
\sl\theaddress\par\medskip \rm Email:\stdspace\tt\theemail\par
\ifx\theurl\relax\else\smallskip \rm URL:\stdspace\tt\theurl\par\fi}}

\def\agtart{
\hoffset 14truemm
\voffset 31truemm
\font\phead=cmsl9 scaled 950
\font\pnum=cmbx10 scaled 913
\font\pfoot=cmsl9 scaled 950
\headline{\vbox to 0pt{\vskip -4.5mm\line{\small\phead\ifnum
\count0=\startpage ISSN numbers are printed here
\hfill {\pnum\folio}\else\ifodd\count0\def\\{ }%
\ifx\theshorttitle\relax\thetitle\else\theshorttitle\fi\hfill{\pnum\folio}
\else\def\\{ and }{\pnum\folio}\hfill\ifx\theshortauthors\relax\theauthors
\else\theshortauthors\fi\fi\fi}\vss}}
\footline{\vbox to 0pt{\vglue 0mm\line{\small\pfoot\ifnum\count0=\startpage
Copyright declaration is printed here\hfill\else
\agt, Volume \thevolumenumber\ (\thevolumeyear)\hfill\fi}\vss}}
\let\maketitle\makeagttitle\let\makeshorttitle\makeagttitle}


\def\ifplaintex{\expandafter\ifx\csname documentclass\endcsname\relax}


\ifplaintex 
\hoffset 14truemm
\voffset 31truemm
\else
\headsep 23pt
\footskip 35pt
\hoffset -4truemm
\voffset 12.5truemm
\fi

\expandafter\ifx\csname beginpicture\endcsname\relax
\expandafter\ifx\csname documentclass\endcsname\relax
\input pictex \else\font\fiverm=cmr5
\input prepictex \input pictex \input postpictex \fi\fi

\def\gt{{\mathsurround=0pt\it $\cal G\mskip-2mu$eometry \&\ 
$\cal T\!\!$opology}}        

\def\gtp{{\mathsurround=0pt\it $\cal G\mskip-2mu$eometry \&\ 
$\cal T\!\!$opology $\cal P\!$ublications}}  


\def\lognumber#1{\def\thelognumber{#1}}
\def\volumenumber#1{\def\thevolumenumber{#1}}
\def\papernumber#1{\def\thepapernumber{#1}}
\def\volumeyear#1{\def\thevolumeyear{#1}}

\def\pagenumbers#1#2{\def\startpage{#1}\def\finishpage{#2}}
\def\published#1{\def\publishdate{#1}}
\def\proposed#1{\def\theproposer{#1}}
\def\seconded#1{\def\theseconders{#1}}
\def\received#1{\def\receiveddate{#1}}
\def\revised#1{\def\reviseddate{#1}}
\def\accepted#1{\def\accepteddate{#1}}
\def\asciititle#1{\def\theasciititle{#1}}

\def\asciiaddress#1{\def\theasciiaddress{#1}}

\long\def\asciiabstract#1{\long\def\theasciiabstract{#1}}
\def\asciikeywords#1{\def\theasciikeywords{#1}}


\let\\\par\let\thelognumber\relax
\let\thevolumenumber\relax\let\thepapernumber\relax
\let\thevolumeyear\relax\let\thesamplenumber\relax\let\startpage\relax
\let\finishpage\relax\let\publishdate\relax\let\receiveddate\relax
\let\reviseddate\relax\let\accepteddate\relax\let\theasciititle\relax
\let\theasciiauthors\relax\let\theasciiaddress\relax
\let\theasciiabstract\relax\let\theasciikeywords\relax
\let\theasciiemail\relax\let\theshortauthors\relax\let\theshorttitle\relax

\long\def\maketitlep{   

\count0=\startpage

\gt\hfill      
\beginpicture
\setcoordinatesystem units <0.33truein, 0.33truein> point at 2.2 0.9
\setplotsymbol ({$\cal G$})
\plotsymbolspacing=9truept
\circulararc 315 degrees from 0 1 center at 0 0
\setplotsymbol ({$\cal T$})
\circulararc 315 degrees from 1 -1 center at 1 0
\endpicture
%
\break
{\small\ifx\thesamplenumber\relax 
Volume \else Sample
\fi\thevolumenumber\ (\thevolumeyear)
\startpage--\finishpage\nl
Published: \publishdate}
\vglue 0.5truein plus 0.4fil minus 0.1truein

{\parskip=0pt\leftskip 0pt plus 1fil\def\\{\par\smallskip}{\ifplaintex\large
\else\Large\fi\bf\thetitle}\par\medskip}   

\vglue 0pt plus 0.1fil 

{\parskip=0pt\leftskip 0pt plus 1fil\def\\{\par}{\sc\theauthors}
\par\medskip}

\vglue 0pt plus 0.1fil 

{\small\parskip=0pt\let\newline\\
{\leftskip 0pt plus 1fil\def\\{\par}{\sl\theaddress}\par}
\expandafter\ifx\theemail\relax    
\relax\else\vglue 5pt plus 0.02fil minus 2pt\def\\{\stdspace{\rm 
and}\stdspace} 
\cl{Email:\stdspace\tt\theemail}\fi
\ifx\theurl\relax                  
\relax\else\vglue 5pt plus 0.02fil minus 2pt\def\\{\stdspace{\rm 
and}\stdspace}
\cl{URL:\stdspace\tt\theurl}\fi\par}

\vglue 7pt plus 0.3fil minus 3pt

{\bf Abstract}
\vglue 5pt plus 0.1fil minus 2pt

\theabstract

\vglue 7pt plus 0.3fil minus 3pt

{\bf AMS Classification numbers}\quad Primary:\quad \theprimaryclass

Secondary:\quad \thesecondaryclass

\vglue 5pt plus 0.3fil minus 2pt

{\bf Keywords}\quad \thekeywords

\vglue 10pt plus 0.5fil minus 5pt

{\small  Proposed: \theproposer\hfill Received: \receiveddate\nl
Seconded: \theseconders\hfill 
\ifx\reviseddate\relax                         
Accepted: \accepteddate                        
\else
Revised: \reviseddate                          
\fi}
\eject
}       

\let\maketitlepage\maketitlep
\let\maketitle\maketitlepage


\font\phead=cmsl9 scaled 950
\font\lhead=cmsl9 scaled 1050
\font\pnum=cmbx10 scaled 913
\font\lnum=cmbx10 
\font\pfoot=cmsl9 scaled 950
\font\lfoot=cmsl9 scaled 1050
\ifplaintex
\headline{\vbox to 0pt{\vskip -4.5mm\line{\small\phead\ifnum
\count0=\startpage ISSN 1364-0380 (on line)
1465-3060 (printed) \hfill {\pnum\folio}\else\ifodd\count0\def\\{ }%
\ifx\theshorttitle\relax\thetitle\else\theshorttitle\fi\hfill{\pnum\folio}
\else\def\\{ and }{\pnum\folio}\hfill\ifx\theshortauthors\relax\theauthors
\else\theshortauthors\fi\fi\fi}\vss}}
\footline{\vbox to 0pt{\vglue 0mm\line{\small\pfoot\ifnum\count0=\startpage
\copyright\ \gtp\hfill\else
\gt, Volume \thevolumenumber\ (\thevolumeyear)\hfill\fi}\vss
}}
\else
\makeatletter
\def\@oddhead{{\small\lhead\ifnum\count0=\startpage ISSN 1364-0380 (on line)
1465-3060 (printed) \hfill {\lnum\number\count0}\else\ifodd\count0
\def\\{ }\ifx\theshorttitle\relax \thetitle \else\theshorttitle\fi\hfill
{\lnum\number\count0}\else\def\\{ and }{\lnum\number\count0}
\hfill\ifx\theshortauthors\relax 
\theauthors\else\theshortauthors\fi\fi\fi}}\def\@evenhead{\@oddhead}
\def\@oddfoot{\small\lfoot\ifnum\count0=\startpage\copyright\ \gtp\hfill\else
\gt, Volume \thevolumenumber\ (\thevolumeyear)\hfill\fi}
\def\@evenfoot{\@oddfoot}
\makeatother
\fi


\newwrite\gtoutfile
\long\gdef\makeheadfile{  
{\def\\{, }\def\s{ }
\immediate\openout\gtoutfile head.xxx
\immediate\write\gtoutfile{To: math@arxiv.org}
\immediate\write\gtoutfile{Subject: put or rep NNNNN:pppp}
\immediate\write\gtoutfile{--text follows this line--}
\immediate\write\gtoutfile{Proxy-for: \ifx\theasciiauthors\relax
\theauthors\else\theasciiauthors\fi\s<\ifx\theasciiemail\relax\theemail\else\theasciiemail\fi>}
\immediate\write\gtoutfile{\noexpand\\}
\immediate\write\gtoutfile{Authors: \ifx\theasciiauthors\relax
\theauthors\else\theasciiauthors\fi}
{\def\\{ }\immediate\write\gtoutfile{Title: \ifx\theasciititle\relax
\thetitle\else\theasciititle\fi}}
\immediate\write\gtoutfile{Subj-class: GT or SG or MG etc}
\immediate\write\gtoutfile{MSC-class: \theprimaryclass\ifx\thesecondaryclass\relax\else, \thesecondaryclass\fi}
\immediate\write\gtoutfile{Journal-ref: Geom. Topol. \thevolumenumber
(\thevolumeyear) \startpage-\finishpage}
\immediate\write\gtoutfile{Comments: Published by Geometry and Topology at}
\immediate\write\gtoutfile{\s\s http://www.maths.warwick.ac.uk/gt/GTVol\thevolumenumber/paper\thepapernumber.abs.html}
\immediate\write\gtoutfile{\noexpand\\}
\immediate\write\gtoutfile{}
\ifx\theasciiabstract\relax
\immediate\write\gtoutfile{\theabstract}\else
\immediate\write\gtoutfile{\theasciiabstract}\fi
\immediate\write\gtoutfile{}
\immediate\write\gtoutfile{\noexpand\\}
\immediate\write\gtoutfile{}
\immediate\closeout\gtoutfile}}  

\def\maketitlepage{\maketitlep\makeheadfile}
\let\maketitle\maketitlepage

\lognumber{238}

\volumenumber{6}
\papernumber{26} 
\volumeyear{2002}
\pagenumbers{889}{904} 
\received{19 February 2002}
\revised{20 December 2002}
\accepted{08 November 2002}
\proposed{Cameron Gordon}
\seconded{Jean-Pierre Otal, Benson Farb}

\published{21 December 2002}

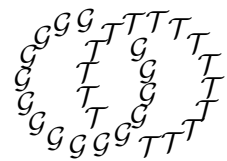
\reflist

\key{1.} {\bf A Basmajian}, {\it Tubular neighbourhoods
of totally geodesic hypersurfaces in hyperbolic manifolds},
Invent. Math. {117} (1994) 207--225

\key{2.} {\bf S Bleiler},  {\bf C Hodgson}, {\it Spherical space forms
and Dehn filling}, Topology {35} (1996) 809--833.

\key{3.} {\bf A Fathi}, {\bf F Laudenbach}, {\bf V Poenaru} 
et al, {\it Travaux de Thurston sur les surfaces},
Asterisque {66--67} (1979)

\key{4.} {\bf S Gallot}, {\bf D Hulin}, {\bf J Lafontaine}, 
{\it Riemannian geometry}, Springer--Verlag (1991)

\key{5.} {\bf J Hass}, {\it Minimal surfaces in manifolds
with $S^1$ actions and the simple loop conjecture for
Seifert fibred spaces}, Proc. Amer. Math. Soc. {99}
(1987) 383--388

\key{6.} {\bf J Hass}, {\bf H Rubinstein}, {\bf S Wang}, {\it
Boundary slopes of immersed surfaces in 3--manifolds}, J.
Differential Geom. {52} (1999) 303--325

\key{7.} {\bf J Hass}, {\bf S Wang}, {\bf Q Zhou}, {\it On
finiteness of the number of boundary slopes of immersed
surfaces in 3--manifolds}, Proc. Amer. Math. Soc. {130}
(2002) 1851--1857

\key{8.} {\bf C Hodgson}, {\bf S Kerckhoff}, {\it
Universal bounds for hyperbolic Dehn surgery}, preprint

\key{9.} {\bf M Lackenby}, {\it Word hyperbolic Dehn surgery},
Invent. Math. {140} (2000) 243--282

\key{10.} {\bf H Masur}, {\it Measured foliations and handlebodies},
Ergod. Th. Dynam. Sys. {6} (1986) 99--116

\key{11.} {\bf J Morgan}, {\it On Thurston's uniformization theorem 
for three-dimensional manifolds}, in The Smith
conjecture (New York, 1979) 37--125, Pure Appl. Math., 112, Academic Press.

\key{12.} {\bf M Scharlemann}, and {\bf Y-Q Wu}, {\it
Hyperbolic manifolds and degenerating handle additions},
J. Austral. Math. Soc. {55} (1993) 72--89

\key{13.} {\bf W Thurston}, {\it The geometry and topology of 
three-manifolds}, Princeton (1980)

\key{14.} {\bf W Thurston}, {\it On the geometry and dynamics
of diffeomorphisms of surfaces}, Bull. Amer. Math. Soc. {19}
(1989) 417--431

\endreflist

\title{Attaching handlebodies to 3--manifolds}
\asciititle{Attaching handlebodies to 3-manifolds}
\author{Marc Lackenby}

\address{Mathematical Institute, Oxford University\\24--29 
St Giles', Oxford OX1 3LB, UK}
\asciiaddress{Mathematical Institute, Oxford University\\24-29 
St Giles', Oxford OX1 3LB, UK}

\email{lackenby@maths.ox.ac.uk}

\abstract
The main theorem of this paper is a generalisation of well known
results about Dehn surgery to the case of attaching handlebodies to a
simple 3--manifold. The existence of a finite set of `exceptional'
curves on the boundary of the 3--manifold is established. Provided none
of these curves is attached to the boundary of a disc in a handlebody,
the resulting manifold is shown to be word hyperbolic and
`hyperbolike'.  We then give constructions of gluing maps satisfying
this condition. These take the form of an arbitrary gluing map
composed with powers of a suitable homeomorphism of the boundary of
the handlebodies.
\endabstract

\asciiabstract{
The main theorem of this paper is a generalisation of well known
results about Dehn surgery to the case of attaching handlebodies to a
simple 3-manifold. The existence of a finite set of `exceptional'
curves on the boundary of the 3-manifold is established. Provided none
of these curves is attached to the boundary of a disc in a handlebody,
the resulting manifold is shown to be word hyperbolic and
`hyperbolike'.  We then give constructions of gluing maps satisfying
this condition. These take the form of an arbitrary gluing map
composed with powers of a suitable homeomorphism of the boundary of
the handlebodies.}

\primaryclass{57N10}\secondaryclass{57N16, 57M50, 20F65}

\asciikeywords{3-manifold, handlebody, word hyperbolic}

\maketitlepage

\section{Introduction}
This paper deals with a generalisation of Dehn surgery. Instead
of attaching solid tori to a 3--manifold with toral boundary components,
we start with a 3--manifold with higher genus boundary components
and glue on handlebodies.
Our aim is to generalise well known surgery results to this setting.
As is customary in surgery theory, we have to assume that the initial
bounded 3--manifold satisfies certain generic topological hypotheses:
it will be {\sl simple}, which means that it is compact, orientable,
irreducible, atoroidal, acylindrical, with incompressible boundary.
Our first and main theorem is the following.
\ppar

\noindent {\bf Theorem 1}\qua {\sl Let $M$ be a simple 3--manifold
with non-empty boundary. Then
there is a finite collection ${\cal C}$ of essential simple closed 
curves on $\partial M$ with the following property. 
Suppose that $H$ is a collection of
handlebodies, and that $\phi \co \partial M \rightarrow \partial H$
is a homeomorphism that sends no curve in ${\cal C}$ to the
boundary of a disc in $H$. Then $M \cup_\phi H$ is
irreducible, atoroidal, word hyperbolic and not Seifert fibred. 
Furthermore, the inclusion map of any component of $H$ into
$M \cup_\phi H$ induces an injection between their
fundamental groups, and hence $\pi_1(M \cup_\phi H)$
is infinite.}\ppar

The set ${\cal C}$ we term the {\sl exceptional curves}. They
may be characterised in terms of the geometry of $M$. For, $M$
has a complete finite volume hyperbolic structure, with
totally geodesic boundary and possibly some cusps. This
may be seen by applying Thurston's geometrisation theorem [11]
to two copies of $M$ glued via the identity map
along their negative Euler characteristic boundary components.
Let $N(\partial M_{\rm cusp})$ be a horoball neighbourhood
of the cusps of $M$, such that, on each component of
$\partial N(\partial M_{\rm cusp})$, the shortest Euclidean
geodesic has length 1. It is shown in [2], for example, 
that $N(\partial M_{\rm cusp})$ is a product neighbourhood of 
the toral ends of $M$.
\ppar

\noindent {\bf Theorem 2}\qua {\sl The set ${\cal C}$ consists of
simple closed geodesics having length at most
$${4\pi \over (1 - 4 /\chi(S))^{1/4} - (1 - 4/\chi(S))^{-1/4} }$$
on boundary components $S$ with genus at least two,
and of closed Euclidean geodesics with length at most
$2\pi$ on $\partial N(\partial M_{\rm cusp})$.
Hence, there is an upper bound on the number of
curves in ${\cal C}$ that depends only on the genus of $\partial M$
and is otherwise independent of $M$.}\ppar

Theorems 1 and 2 can be viewed as a generalisation of Thurston's hyperbolic
Dehn surgery theorem [13] and its extension by Hodgson and
Kerckhoff [8]. Thurston established the existence of
a finite set ${\cal C}$ of slopes on a torally bounded hyperbolic
3--manifold, such that, provided these slopes are avoided when surgery
is performed, the result is a 3--manifold with a hyperbolic structure.
Hodgson and Kerckhoff provided a universal upper bound,
independent of the 3--manifold, on the number of curves in 
${\cal C}$ on each torus boundary component.

Theorems 1 and 2 should be compared with the main theorem of [12],
due to Scharlemann and Wu. They considered the attachment of
a single 2--handle to $M$. One is only allowed to attach it along
a certain type of curve, known as a `basic' curve, but this
is not a serious restriction. They prove that the resulting
manifold is hyperbolic, provided one avoids a finite set of
curves on $\partial M$, up to isotopy. Attaching a handlebody
can be performed by first gluing on 2--handles along basic
curves and then Dehn filling. Thus, their result implies that,
in a certain sense, most ways of attaching a handlebody
give a hyperbolic manifold. The limitation of Scharlemann and Wu's
procedure is that the 2--handles must be attached in sequence,
and then be followed by the surgeries. Thus, the curves that
the later 2--handles and surgeries must avoid depend on
where the earlier 2--handles are attached. In Theorem 1,
we attach the entire handlebody in a single step.
This allows us to identify the exceptional curves at
the outset.

Of course, Theorem 1 raises the problem of finding homeomorphisms
$\phi$ with the required property.  The approach we consider is to
start with an arbitrary $\phi$, and then modify this by applying
powers of a homeomorphism $f  \co \partial H \rightarrow \partial H$.
If $f$ extends to a homeomorphism of $H$, then this will not change the
resulting manifold. More generally, if some power of $f$ extends to
a homeomorphism, then only finitely many manifolds are created.
We would like to avoid this situation.
In fact, one must rule out a yet more general possibility:
no power $f^n$ of $f$ (where $n \not= 0$)
may {\sl partially extend to} $H$. This
means that there is a compression body $R$ (other than
a product) embedded in
$H$ with positive boundary $\partial H$, such that $f^n$
extends to a homeomorphism of $R$. For, if $\phi$ were to
map a curve $C$ in ${\cal C}$ to the boundary of a disc in $R$,
then $f^n \phi(C)$ would bound a disc in $R$ for infinitely many
$n$. However, if this condition is met, we shall show
that this method does create infinitely many manifolds that satisfy
the conclusions of Theorem 1.

\ppar
\noindent {\bf Theorem 3}\qua {\sl Let $M$ be a simple 3--manifold,
let $H$ be a collection of handlebodies, and let $\phi \co
\partial M \rightarrow \partial H$ be a homeomorphism. Let
$f \co \partial H \rightarrow \partial H$ be a
homeomorphism, no power of which partially extends to
$H$. Then, for infinitely many integers $n$, $M \cup_{f^n \circ \phi} H$
satisfies the conclusions of Theorem 1.}\ppar

An alternative method of finding suitable homeomorphisms $f$
is to use the theory of pseudo-Anosov maps ([14], [3]). We can regard
simple closed curves on $\partial H$, with counting measure, as elements
of the space $PL(\partial H)$ of projective measured
laminations on $\partial H$. Let $B(\partial H)$ be
the closure in $PL(\partial H)$ of the set of curves
on $\partial H$ that bound discs in $H$. It is a theorem of Masur [10]
that $B(\partial H)$ is nowhere dense in $PL(\partial H)$.
Hence that the stable and unstable laminations
of a `generic' pseudo-Anosov homeomorphism $f \co
\partial H \rightarrow \partial H$ will not lie
in $B(\partial H)$.

\ppar
\noindent {\bf Theorem 4}\qua {\sl Let $M$ be a simple 3--manifold,
let $H$ be a collection of handlebodies, and let $\phi \co
\partial M \rightarrow \partial H$ be a homeomorphism.
Let $f \co \partial H \rightarrow \partial H$ be a
pseudo-Anosov homeomorphism whose stable and unstable
laminations do not lie in $B(\partial H)$. Then for
all but finitely many integers $n$, $M \cup_{f^n \circ \phi} H$
satisfies the conclusions of Theorem 1.}\ppar

The proof of Theorems 1 and 2 follows two papers: [9],
which established the word hyperbolicity of certain surgered
manifolds, and [7] by Hass, Wang and Zhou, which examined 
boundary slopes of immersed essential surfaces in hyperbolic 
3--manifolds with totally geodesic
boundary. Theorem 3 uses some new techniques, involving
arrangements of discs in a handlebody and a delicate counting
argument. Theorem 4 is
an elementary application of the theory of pseudo-Anosov
automorphisms.

This paper suggests many interesting areas for further research.
Firstly, can the results be upgraded to deduce the existence of
a metric that is hyperbolic or just negatively curved? Secondly, 
can Theorem 3 be strengthened
so that the conclusion holds for all but finitely many $n$?
Thirdly, can the techniques of this paper be generalised to
analyse Heegaard splittings?

\section{The main theorem}

In this section, we will prove Theorems 1 and 2. Suppose therefore that
$M$ is a simple 3--manifold with non-empty boundary. Denote its
toral boundary components by $\partial M_{\rm cusp}$ and the
remaining boundary components by $\partial M_{\rm geo}$. We work with
the complete finite volume hyperbolic structure on
$M - \partial M_{\rm cusp}$ in which $\partial M_{\rm geo}$
is totally geodesic. Let ${\cal C}$ be the set of simple closed
curves on $\partial M$ as described in Theorem 2.
The fact that there is an upper bound on the number of
curves in ${\cal C}$ that depends only on the genus of
$\partial M$ is well known. A proof is given in [7] for
example. Suppose that $\phi \co \partial M \rightarrow \partial H$ is
a homeomorphism that sends no curve in ${\cal C}$ to the boundary
of a disc in $H$.

Our first step is to show that the inclusion map of
any component of $H$ into $M \cup_\phi H$ induces an
injection between their fundamental groups. Consider
an essential loop $L$ in $H$, and suppose that it is homotopically
trivial in $M \cup_\phi H$. There is then a map $f \co D \rightarrow
M \cup_\phi H$, where $D$ is a disc, such that $f|_{\partial D}$ winds
once around $L$. Homotope $f$ a little, so that $f^{-1}(\partial H)$ is
a collection of simple closed curves in the interior of $D$, and so
that $f$ is transverse to $\partial H$ near these curves. We
suppose that $L$ and $f$ have been chosen so that the number of
these curves has been minimised.

\proclaim{Claim 1} $f^{-1}(H)$ is a collection of discs
and a collar on $\partial D$.\endproc

If not, pick a curve $L'$ of $f^{-1}(\partial H)$ that is
innermost in $D$ among curves that do not bound discs of $f^{-1}(H)$
and that are not parallel in $f^{-1}(H)$ to $\partial D$.
It bounds a disc $D'$. If $L'$ is homotopically trivial in $H$,
we may modify $f$ in $D'$ so that it is mapped entirely to $H$,
thus reducing $|f^{-1}(\partial H)|$,
which contradicts the minimality of $|f^{-1}(\partial H)|$. If
$L'$ is homotopically non-trivial in $H$, we may work with $D'$ instead
of $D$. Again the minimality assumption is violated.

\proclaim{Claim 2} The surface $F = f^{-1}(M)$ is homotopically 
boundary-incompressible in $M$.\endproc 

Recall from [9] that this means that no properly embedded essential arc in $F$
can be homotoped in $M$, keeping its endpoints fixed, to an arc in $\partial M$.
For, if there were such an arc $\alpha$, we could perform a homotopy 
to $f$, taking a regular neighbourhood of $\alpha$ into $H$. There are two
cases to consider: when the endpoints of $\alpha$ lie in distinct
boundary components of $F$, and when they lie in the
same component. In the first case, the result is to reduce
$|f^{-1}(\partial H)|$, contradicting the minimality assumption.

In the second case, there are two subcases: either $\partial \alpha$
lies in the collar on $\partial D$ or it does not. Suppose first that $\partial \alpha$
misses the collar. The arc separates $F$ into two components,
as $F$ is planar. One of these, $F'$ say, does not intersect the
collar on $\partial D$. Let $L'$ be the boundary component
of $F'$ that runs along the neighbourhood of $\alpha$. Push
$L'$ a little into $H$. Since $F'$ has fewer boundary components
than $F$ and $L'$ is homotopically trivial in $M \cup_\phi H$, the
minimality assumption implies that $L'$ is homotopically trivial in
$H$. We may therefore remove $F'$ from $F$ and replace it with
a disc that maps to $H$. This reduces $|\partial F|$, which
is a contradiction.

Suppose now that both endpoints of $\alpha$ lie in the
collar on $\partial D$.
We consider the two halves of $F$ cut open along $\alpha$. By the
minimality assumption, the boundary curves of each are homotopically
trivial in $H$. But then $L$ is trivial in $H$, which is a 
contradiction.

\proclaim{Claim 3} $F$ is homotopically
incompressible in $M$.\endproc

This means that no homotopically non-trivial simple closed
curve in $F$ maps to a homotopically trivial curve in $M$.
The argument is similar to that of Claims 1 and 2,
but simpler, and so is omitted.

We now follow the argument of [7]. 
Let $N(\partial M_{\rm cusp})$ be a horoball neighbourhood
of the cusps of $\partial M$, such that on each component
of $\partial N(\partial M_{\rm cusp})$, the shortest
geodesic has length 1. Let $N(\partial M_{\rm geo})$
be the set of points at a distance at most $U$
from $\partial M_{\rm geo}$. By a theorem of Basmajian [1], if we take 
$U$ to be 
$${1 \over 4} \log \left ( 1- {4 \over \chi(S)} \right )
= \sinh^{-1} \left ( {(1 - 4 / \chi(S))^{1/4} - (1 - 4/\chi(S))^{-1/4} \over 2 }
\right),$$
then $N(\partial M_{\rm cusp}) \cup N(\partial M_{\rm geo})$ will be a collar
on $\partial M$. Denote $N(\partial M_{\rm cusp})
\cup N(\partial M_{\rm geo})$ by $N(\partial M)$.

\proclaim{Claim 4.} There is a least area surface in the
homotopy class of $f \co (F, \partial F) \rightarrow
(M, \partial M)$. This is an immersion.\endproc

If $F$ were closed and $M$ had no cusps, this would be Lemma 2
of [5], as $F$ is homotopically incompressible, by Claim 3.
This was extended to the case where $F$ has boundary
and $M$ has cusps in Theorem 4.4 of [6]. In [6],
though, it was assumed that $f_\ast \co
\pi_1(F) \rightarrow \pi_1(M)$ and $f_\ast \co
\pi_1(F,\partial F) \rightarrow \pi_1(M,\partial M)$
are injective. However, as explained in [5], these
hypotheses can be weakened to the assumption that
$F$ is homotopically incompressible and homotopically
boundary-incompressible.

We perform the homotopy in Claim 4. Let 
$\partial F_{\rm cusp}$ and $\partial 
F_{\rm geo}$ be the boundary components of $F$ that map
to $\partial M_{\rm cusp}$ and $\partial M_{\rm geo}$
respectively. Then $F - \partial F_{\rm cusp}$ 
inherits a Riemannian metric with geodesic boundary.
Its sectional curvature is the sum of the
sectional curvature of $M$ and the product of its
principal curvatures (see Theorem 5.5 of [4], for example).
The principal curvatures sum to zero, since $F$
is minimal, and hence their product is non-positive.
Therefore, the sectional curvature of $F$ is at
most $-1$.

We abuse notation slightly by denoting the
component of $\partial F$ parallel to $L$ by $L$.
Let $N(\partial F_{\rm cusp})$
be those components of the inverse image of 
$N(\partial M_{\rm cusp})$ in $F$ that contain a component
of $\partial F_{\rm cusp}$.  Let $N(\partial F_{\rm geo})$ be the set of points in $F$
with distance at most $U$ from $\partial F_{\rm geo}$ in its intrinsic path metric. 
Denote $N(\partial F_{\rm cusp}) \cup N(\partial F_{\rm geo})$ by
$N(\partial F)$.

\proclaim{Claim 5} $N(\partial F)$ is a collar on $\partial F$.\endproc

We start by showing that each component of $N(\partial F_{\rm geo})$
is a collar. To see this, increase $U$ from zero to its final
value. Near zero, $N(\partial F_{\rm geo})$ is clearly a
collar. But suppose that as it expands, there is some
point at which a self-tangency is created. Then there
is a geodesic arc properly embedded in $N(\partial F_{\rm geo})$
with endpoints perpendicular to $\partial F_{\rm geo}$.
This geodesic is the concatenation of two geodesic
arcs which run from $\partial F_{\rm geo}$ to the
point of self-tangency.
Since $F$ is negatively curved, this arc is essential in $F$.
But it lies within $N(\partial M_{\rm geo})$, which is
a collar, and so can be homotoped into $\partial M$, keeping
its endpoints fixed. This contradicts the fact that $F$
is homotopically boundary-incompressible.

It is also clear that each component of
$N(\partial F_{\rm cusp})$ is a collar. For,
otherwise, it contains a properly embedded arc that is essential
in $N(\partial F_{\rm cusp})$ and that has endpoints in
$\partial F_{\rm cusp}$. The fact that $N(\partial M_{\rm cusp})$
is a collar on $\partial M_{\rm cusp}$ implies that this arc can be
homotoped into $\partial M$. The arc is therefore
inessential in $F$. Hence, some component of $\partial N(\partial F_{\rm
cusp})$ bounds a disc in $F$. Consider the lift of
this disc to ${\Bbb H}^3$, which contains the universal cover of $M$.
The boundary of the disc lies on a component of the inverse image
of $\partial N(\partial M_{\rm cusp})$, which is a
horotorus. There is a natural projection from
${\Bbb H}^3$ onto this horotorus, and this reduces
the area of the disc. This contradicts the assumption
that $F$ is least area.

Finally, $N(\partial F_{\rm cusp})$ and $N(\partial F_{\rm geo})$
must be disjoint. For the former lies in $N(\partial M_{\rm cusp})$,
whereas the latter lies in $N(\partial M_{\rm geo})$, and these
are disjoint. This proves the claim.

\proclaim{Claim 6} ${\rm Area}(N(\partial F_{\rm geo}))
\geq (\sinh U) \ ({\rm Length}(\partial F_{\rm geo}))$.\endproc

Choosing orthogonal co-ordinates for the surface $N(\partial F_{\rm geo})$, the
metric is given by
$$ds^2 = du^2 + J^2(u,v) dv^2,$$
where $J(u,v) > 0$ and $J(0,v) = 1$. The curves
where $v$ is constant are geodesics perpendicular
to the boundary, and the curve $u = 0$
lies in the boundary. By the calculations in [7], 
$${\partial J \over \partial u} \geq 0, \qquad 
{\partial^2 J \over \partial u^2} \geq J,$$
and hence $J(u,v) \geq \cosh(u)$. Therefore,
$$\eqalign{
{\rm Area}(N(\partial F_{\rm geo}))
&= \int_{v=0}^{{\rm Length}(\partial F_{\rm geo})}
\int_{u=0}^U J(u,v) \ du \ dv \cr
&\geq \int_{v=0}^{{\rm Length}(\partial F_{\rm geo})}
\int_{u=0}^U \cosh u \ du \ dv \cr
&= (\sinh U) \ ({\rm Length}(\partial F_{\rm geo})).}$$

Give $\partial H$ the Riemannian metric on $\partial N(\partial M_{\rm
cusp}) \cup \partial M_{\rm geo}$.
Each component of $\partial F$, except the one parallel to
$L$, is mapped to the boundary
of a disc in $H$, by Claim 1. It is essential in $\partial H$
by Claim 3. So, by assumption, it has length more than the
bound in Theorem 2 if it is simple. The non-simple case follows 
from the following claim.

\proclaim{Claim 7} Each shortest essential closed curve in $\partial H$
that is homotopically trivial in $H$ is simple.\endproc

Let $\tilde H$ be the universal cover
of $H$. Then $\partial \tilde H$ is the cover of $\partial H$
corresponding to the kernel of $\pi_1(\partial H) \rightarrow
\pi_1(H)$. In particular, it is a regular cover. The closed curves
in $\partial \tilde H$ are exactly the lifts of curves in
$\partial H$ that are homotopically trivial in $H$. The shortest 
of these, $C$ say, that is essential in $\partial \tilde H$
is a geodesic. It must be simple, for otherwise, it can be
modified at a singular point to reduce its length. 

We claim that $C$
must be disjoint from its covering translates. For if $C$
and some translate $C'$ were to intersect transversely, they would do so
at least twice, since both are homologically trivial in $\tilde H$.
Pick two points of intersection, which divide $C$ and $C'$ each into
two arcs. Glue the shorter of the two arcs in $C$ (or either of these arcs
if they have the same length) to the shorter of the two in $C'$,
and then smooth off to form a shorter closed curve
in $\partial \tilde H$. If this is essential in $\partial \tilde H$,
we have a contradiction. If not, then we may homotope the
longer of its two arcs onto the other, and then smooth off.
This reduces the length of $C$ or $C'$, which again is a
contradiction.

If $C$ and $C'$ were to intersect non-transversely, then they would
be equal. Hence, $C$ would project to a multiple of a simple
closed curve $C''$ in $\partial H$. If we apply the Loop Theorem
to a regular neighbourhood of $C''$, we see that $C''$
must also bound a disc in $H$. But this is a shorter curve
than $C$, which is a contradiction.

Hence, $C$ projects homeomorphically to a simple closed
curve in $\partial H$ which is shortest among essential
curves that are trivial in $H$. This proves the claim.

Putting together Claims 6 and 7, we obtain
$${\rm Area}(N(\partial F_{\rm geo})) \geq 2\pi \, |\partial F_{\rm
geo} - L|.$$

Let ${\rm Length}(\partial F_{\rm cusp})$ be the length
of the curves $\partial N(\partial F_{\rm cusp})$ 
on $\partial N(\partial M_{\rm cusp})$.
By Claim 5, each of these curves has the same slope as the
corresponding component of $\partial F$. Hence, each
(except possibly $L$) has length more than $2 \pi$, by assumption. Thus,
using a well known argument, which can be found in 
the proof of Theorem 4.3 of [6], for example, we obtain the following.

\proclaim{Claim 8} ${\rm Area}(N(\partial F_{\rm cusp}))
\geq 2 \pi \, |\partial F_{\rm cusp} - L|$.\endproc

Since $F$ is a minimal surface in a hyperbolic manifold, 
its sectional curvature $K$ is at
most $-1$. So, applying Gauss--Bonnet to $F$:
$$\eqalign{& 2\pi (2 - |\partial F|) =
2 \pi \chi(F) =
\int_{F} K \ dA
\leq -{\rm Area}(F) \cr
&\qquad \leq -{\rm Area}(N(\partial F_{\rm cusp}))
-{\rm Area}(N(\partial F_{\rm geo}) )
\leq -2\pi (|\partial F| - 1),}
$$
which is a contradiction. We therefore deduce that
the inclusion of each component of $H$ into $M \cup_\phi H$
is $\pi_1$--injective.

A very similar argument can be used to show that $M \cup_\phi H$
is irreducible. For otherwise, $M$ would contain a properly
embedded incompressible boundary-incompressible planar surface, each
boundary component of which would extend to a disc in $H$.
(To make this argument work, we need the fact that
the inclusion of $H$ into $M \cup_\phi H$ is $\pi_1$--injective.)
We can then apply the above analysis to this surface.

It remains to show that $M \cup_\phi H$ is word hyperbolic.
For this implies that $\pi_1(M \cup_\phi H)$ contains no
rank two abelian subgroup, and hence that $M \cup_\phi H$
is atoroidal. The fact that $\pi_1(M \cup_\phi H)$
is infinite then gives that $M \cup_\phi H$ is not Seifert fibred.

Pick a Riemannian metric $g$ on $M \cup_\phi H$, in which
$H$ is an $\epsilon$--neighbourhood of a graph, for small $\epsilon >0$.
We wish to establish a linear isoperimetric inequality for $g$.
So, consider a loop $L$ that is homotopically trivial,
and which therefore forms the boundary of a mapped-in disc
$f \co D \rightarrow M \cup_\phi H$. Using a small homotopy,
we may move $L$ off $H$. We can ensure that this changes
the length of $L$ by a factor that is bounded independently
of $L$. We can also ensure that the area of the annulus
realizing the homotopy is similarly bounded.

By the argument of Claims 1, 2 and 3 and using the fact that
the inclusion of $H$ into $M \cup_\phi H$ is $\pi_1$--injective, we may
homotope $f$, keeping it fixed on $L$, so that 
afterwards $F = f^{-1}(M)$ 
is homotopically incompressible and homotopically boundary-incompressible 
in $M$, and so that $f^{-1}(H)$ is a collection of discs.

Let $h$ be the hyperbolic metric on $M - \partial M_{\rm cusp}$. 
Since $M$ is compact, the Riemannian manifolds
$(M - {\rm int}(N(\partial M_{\rm cusp})), h)$ 
and $(M,g)$ are bi-Lipschitz equivalent, with constant $c_1$, say. 
Also, there is a map $(N(\partial M_{\rm cusp}), h) \rightarrow
(H,g)$, collapsing the cusps to solid tori, that
increases areas by at most a factor $c_1^2$, say.
Initially, we will measure lengths
and areas in the $h$ metric. Note that, by construction,
$L$ is disjoint from $N(\partial M_{\rm cusp})$.

We may assume that $L$ is homotopically non-trivial in $M$.
For otherwise it is trivial to construct a disc in $M$ bounded by $L$
with area linearly bounded by the length of $L$. 
(Indeed, the fact that closed curves in hyperbolic
space have this property is the motivation behind
Gromov's theory of hyperbolic groups.) Also,
we will, for the moment, assume that $L$ is not homotopic
to a curve in $\partial M_{\rm cusp}$.
Hence it has a geodesic representative $\overline L$. 
We may realize the free homotopy between
$L$ and $\overline L$ with a mapped-in annulus $A$
having area at most $c_2 {\rm Length}(L)$, for some
constant $c_2$ depending only on $h$. We may assume
that $A$ has least possible area.

\proclaim{Claim 9} Each component of $A \cap N(\partial M_{\rm
cusp})$ is either a disc disjoint from $\overline L$, or a disc
containing a single component of $\overline L \cap N(\partial M_{\rm
cusp})$.\endproc

Note that $A \cap N(\partial M_{\rm cusp})$ is disjoint
from $L$, and does not contain a core curve of $A$.
Suppose that the claim is not true. Then we can find a component of
$A \cap N(\partial M_{\rm cusp})$ which is not a disc
or which intersects $\overline L$ more than once.
In the former case, $A \cap N(\partial M_{\rm cusp})$ has a boundary component
that is a simple closed curve bounding a disc in $A$
with interior disjoint from $N(\partial M_{\rm cusp})$. However, 
we could then homotope this disc into $N(\partial M_{\rm cusp})$
to reduce the area of $A$, which is a contradiction.
Thus, each component of $A \cap N(\partial M_{\rm cusp})$
is a disc. If one component intersects $\overline L$
in more than one arc, then the sub-arc of
$\overline L$ between these two arcs does not lie wholly in
$N(\partial M_{\rm cusp})$, but can be homotoped
into $N(\partial M_{\rm cusp})$. This contradicts the
fact that $\overline L$ is a geodesic. This proves
the claim.

Homotope $F$ (minus those components that map to
$\partial M_{\rm cusp}$) to a least area surface
in $M - \partial M_{\rm cusp}$, with one
boundary component mapping to $\overline L$,
and the remainder mapping to $\partial M$.
Divide the boundary components of $F$ into three subsets
$\overline L$, $\partial F_{\rm cusp}$ and $\partial F_{\rm geo}$.
Note that $F \cup A$ is now a planar surface with one
boundary component mapping to $L$ and the remaining components
mapping to curves in $\partial M$ that bound discs
in $H$. This surface will extend to a disc in
$M \cup_\phi H$ with area linearly bounded by
the length of $L$, establishing the required
linear isoperimetric inequality.

We define $N(\partial F_{\rm cusp})$ and $N(\partial F_{\rm geo})$
as before. Claim 5 now reads that\break $N(\partial F_{\rm cusp}) \cup
N(\partial F_{\rm geo})$ is a collar on $\partial F_{\rm cusp}
\cup \partial F_{\rm geo}$, which may have non-empty intersection
with $\overline L$. However, the 
calculations in Claims 6 and 8 now need modification.

Further divide $\partial F_{\rm geo}$
into two subsets: $\partial F_{\rm thin}$ and $\partial F_{\rm thick}$.
The former is the set of points $x$ in $\partial F_{\rm geo}$ such that,
when a perpendicular geodesic is emitted from $x$ in $F$, it
meets $\overline L$ within a distance $U$. Let $\partial F_{\rm thick}$
be the remainder of $\partial F_{\rm geo}$. Then the argument of Claim 6
gives that
$${\rm Area}(N(\partial F_{\rm geo})) \geq (\sinh U) \ ({\rm Length}(\partial F_
{\rm thick})).$$
By Claim 9,
$N(\partial F_{\rm cusp}) \cup (A \cap N(\partial M_{\rm cusp}))$ 
is a collar on $\partial F_{\rm cusp}$ and possibly
some discs. Using these collars, we may associate
to each component of $\partial F_{\rm cusp}$
a curve on $\partial N(\partial M_{\rm cusp})$,
namely the relevant boundary component of 
$N(\partial F_{\rm cusp}) \cup (A \cap N(\partial M_{\rm cusp}))$.
We define the length of $\partial F_{\rm cusp}$
to be the length of these curves. Then, Claim 8, 
suitably modified, gives that
$${\rm Area}(N(\partial F_{\rm cusp}) \cup A) \geq 
{\rm Length}(\partial F_{\rm cusp}),$$
and hence that
$${\rm Area}(N(\partial F_{\rm cusp})) \geq
{\rm Length}(\partial F_{\rm cusp}) - c_2 \ {\rm Length}(L).$$
Therefore,
$$\eqalignno{
{\rm Area}(F) &\geq {\rm Area}(N(\partial F_{\rm geo}))
+ {\rm Area}(N(\partial F_{\rm cusp})) \cr
&\geq 
(\sinh U) \ ({\rm Length}(\partial F_{\rm thick}))\cr
&\qquad + {\rm Length}(\partial F_{\rm cusp}) - c_2 \ {\rm Length}(L). &(1)}$$
Also,
$${\rm Length}(L) \geq {\rm Length}(\overline L) \geq
{\rm Length}(\overline L \cap N(\partial M_{\rm geo}))
\geq {\rm Length}(\partial F_{\rm thin}).\eqno{(2)}$$
Let $c_3$ (respectively, $c_4$) denote the minimal length
of a geodesic on $\partial N(\partial M_{\rm cusp})$
(respectively, $\partial M_{\rm geo}$) with length more
than $2\pi$ (respectively, $2 \pi /\sinh(U)$), and let
$c_5$ be $\min \{ c_3, c_4 \sinh(U) \}$, which is
more than $2\pi$. Then, since each curve of $\partial F_{\rm cusp}$
and $\partial F_{\rm geo}$ has length more than the bound
in Theorem 2, we know that
$$\eqalignno{
& {\rm Length}(\partial F_{\rm cusp}) +
(\sinh U)({\rm Length}(\partial F_{\rm geo})) \cr
&\qquad \geq c_5(|\partial F| - 1) >
-c_5 \ \chi(F) \geq c_5 \ {\rm Area}(F) / 2\pi. &(3)
\cr}$$
The final inequality is an application of Gauss--Bonnet.
So, adding $(2 \pi / c_5)$ times (3),
and $(\sinh U)$ times (2), to (1), and cancelling the area term,
we get:
$$\eqalign{
&(c_2 + \sinh U){\rm Length}(L) \cr
&\qquad \geq (1-2 \pi / c_5)
\left ( 
{\rm Length}(\partial F_{\rm cusp})
+(\sinh U)({\rm Length}(\partial F_{\rm geo})) \right ),}$$
which we summarise as
$${\rm Length}(L)
\geq c_6 \ {\rm Length}(\partial F_{\rm geo} \cup \partial F_{\rm
cusp}),$$
for some positive constant $c_6$ independent of $L$.
So, by equation (3),
$${\rm Length}(L) \geq
c_7 \ {\rm Length}(\partial F_{\rm geo} \cup \partial F_{\rm cusp})
+ c_8 \ {\rm Area}(F),\eqno{(4)}$$
where $c_7$ and $c_8$ are positive constants independent of $L$.

Recall that we assumed earlier that $L$ was not homotopic in
$M$ to a curve in $\partial M_{\rm cusp}$. We may now drop
this assumption, since in this case, we can easily find a
surface $F$ satisfying the above inequality.
One method is as follows. Let $L'$ be the curve in
$\partial N(\partial M_{\rm cusp})$ that is homotopic
to $L$, and that is a geodesic in the Euclidean Riemannian metric
on $\partial N(\partial M_{\rm cusp})$. Then $L'$
has length less than that of $L$. The homotopy between
$L$ and $L'$ is realized by an annulus $A'$, say.
Identify $A'$ with $S^1 \times I$, and homotope
each arc $\{ \ast \} \times I$, keeping its
endpoints fixed, to a geodesic. It is not hard to
calculate that the area of $A'$ after this
homotopy is at most $({\rm Length}(L) +
{\rm Length}(L')) \leq 2 \ {\rm Length}(L)$.
Now let $F$ be the union of $A$ and the vertical
annulus above $L'$ in $N(\partial M_{\rm cusp})$.
Then, the area of $F$ is at most $3 \ {\rm Length}(L)$.
Also, the length of $\partial F_{\rm geo} \cup \partial F_{\rm cusp}$ is 
the length of $L'$, which is less than that of $L$. 
This verifies (4) in this case.

Changing (4) to the metric $g$, we get
$${\rm Length}(L,g) \geq (c_7/c_1^2) {\rm Length}(\partial F_{\rm geo} \cup \partial
F_{\rm cusp},g) + (c_8/c_1^3) {\rm Area}(F,g).$$
Now, each component of $H$ has free fundamental group, and hence is word hyperbolic.
So, the curves $\partial F_{\rm geo} \cup \partial F_{\rm cusp}$ bound discs in $H$
with area at most $c_9 \, {\rm Length}(
\partial F_{\rm geo} \cup \partial F_{\rm cusp},g)$,
for some constant $c_9$.
Attaching these discs to $F$, and then attaching the annulus $A$
between $L$ and $\overline L$, we create a disc bounded by $L$
with area at most $c_{10} \, {\rm Length}(L,g)$, for some constant
$c_{10}$ independent of $L$. This establishes the required linear isoperimetric
inequality, and hence the proof of Theorems 1 and 2 is complete.

\section{Modifying the gluing map by powers of a homeomorphism}

In this section, we prove Theorem 3. Therefore let $\phi \co
\partial M \rightarrow \partial H$ be some homeomorphism, and let
$f \co \partial H \rightarrow \partial H$ be a
homeomorphism, no power of which partially extends to $H$.
Let ${\cal C} = \{ C_1, \dots, C_{|{\cal C}|} \}$ be the exceptional
curves on $\partial M$. Suppose that Theorem 3 does not hold.
Then, there is some finite subset ${\cal S}$ of ${\Bbb Z}$,
and a function $i \co {\Bbb Z} - {\cal S}
\rightarrow \{ 1, \dots, |{\cal C}| \}$ such that for
each $n \in {\Bbb Z} - {\cal S}$,
$f^n \phi(C_{i(n)})$ bounds a disc $D_n$ in $H$.

We wish to analyse these discs $D_n$ in $H$, and
we therefore introduce a definition. A {\sl disc arrangement} is a
collection of properly embedded discs in general position
in a 3--manifold. Two discs arrangements are {\sl equivalent}
if there is a homeomorphism between their regular neighbourhoods
taking one set of discs to the other. We term this a
{\sl relative equivalence} if the homeomorphism is the
identity on the boundary of the discs. A {\sl double curve}
is a component of intersection between two of the discs. 
The {\sl boundary curves}
of a disc arrangement are simply the boundaries of the
discs.

\ppar{\bf Lemma 1}\qua {\sl Fix a finite collection of embedded (not
necessarily disjoint) simple
closed curves in general position in a surface $S$. Suppose that
these form the boundary curves of a disc arrangement in some 
irreducible 3--manifold bounded by $S$.
Then we may isotope these discs keeping their boundaries fixed, 
so that afterwards, they belong to
one of only finitely many relative equivalence classes of
disc arrangements, and so that a regular neighbourhood of
these discs and $S$ is homeomorphic to a punctured compression body 
with positive boundary $S$.}

\prf This is by induction
on partition number, which we define to be the smallest number of
subsets into which we may partition the collection
of discs, such that if two discs belong to the same
subset, they are disjoint. When the partition number is one, the
lemma is trivial, and the induction is started.
We now prove the inductive hypothesis. Let ${\cal D}$ denote
the collection of discs, and let ${\cal D}_0$ be one of the subsets 
in a partition of ${\cal D}$ that realizes the partition number.

First perform an isotopy to the discs not in ${\cal D}_0$
to remove any simple closed double curves that lie in
$\bigcup {\cal D}_0$. This is achieved by an elementary
innermost curve argument.

We now perform isotopies on the discs not in ${\cal D}_0$,
supported in a regular neighbourhood of $\bigcup {\cal D}_0$,
so that afterwards any two double curves in $\bigcup {\cal D}_0$
intersect in at most one point. For if not, two double curves
form a bigon in some disc in ${\cal D}_0$. We may
remove an innermost bigon by an isotopy. Note that this
introduces no simple closed double curves in $\bigcup {\cal D}_0$.
After this isotopy, there are only finitely many possibilities
for the double curves in $\bigcup {\cal D}_0$, up to ambient
isotopy fixed on the boundary of $\bigcup {\cal D}_0$.

Throughout these isotopies, discs which were initially disjoint
have remained so. Thus, the partition number has not increased.
Cut the ambient manifold along $\bigcup {\cal D}_0$
to a give a 3--manifold with boundary $S'$, say. The remaining
discs of ${\cal D}$ are cut up to form a disc arrangement
${\cal D}'$ in the cut-open 3--manifold. Its boundary curves arise from the boundary 
curves of ${\cal D} - {\cal D}_0$ and from the double curves
in $\bigcup {\cal D}_0$. Thus, there are only finitely many
possibilities for these boundary curves.

Now, ${\cal D}'$ has lower partition number than the
original ${\cal D}$. Hence, inductively, ${\cal D}'$ may be isotoped,
keeping its boundary fixed, so that afterwards, it belongs to
one of only finitely many relative equivalence classes of
disc arrangements. Thus, the same is true of ${\cal D}$. 
Also, a regular neighbourhood of $S' \cup \bigcup
{\cal D'}$ is a punctured compression body with
positive boundary $S'$, and so a regular neighbourhood
of $S \cup \bigcup {\cal D}$ is as required. \endprf

We will always take any given finite subcollection of the discs $D_n$
to be one of the disc arrangements as in Lemma 1. 

We define the {\sl complexity} of a closed orientable
surface to be twice its genus, minus the number of
its non-spherical components. If it is a union of
2--spheres or it is empty, we take its complexity to
be zero. Thus, complexity is a non-negative integer, with
the property that if the surface is
compressed, its complexity goes down. So, if one
compression body $R'$ is embedded within another, $R$,
and they have the same positive boundary, then the
complexity of $\partial_- R$ is at most the complexity
of $\partial_- R'$. This is an equality if and only
if $R - R'$ is a regular neighbourhood of a subsurface
of $\partial_- R$.

Let $s = \max \{ n : n \in {\cal S} \}$. For integers
$s < a \leq b$, we define $H(a,b)$ to be the submanifold of 
$H$ formed by taking a regular neighbourhood of 
$\partial H \cup \bigcup_{a \leq n \leq b} D_n$,
and then filling in any 2--sphere boundary components
with 3--balls. This is a compression body with positive
boundary $\partial H$. Let $h(a,b)$ be the complexity
of its negative boundary.
Note that if $a \leq a'$ and $b \leq b'$,
then $h(a,b)$ and $h(a',b')$ are at least $h(a,b')$.

We term two pairs of integers $(a,b)$ and $(a',b')$
{\sl comparable} if $f^{a'-a}$ takes the curves
$\partial D_a, \dots, \partial D_b$
to the curves $\partial D_{a'}, \dots, \partial D_{b'}$,
and this extends to an equivalence of disc arrangements.
Note that in this case, $f^{a'-a}$ extends to a
homeomorphism of $H(a,b)$ to $H(a',b')$.

We claim that we can find comparable pairs of
integers $(a,b)$ and $(a',b')$, where $a < a'$, such that
$h(a,b) = h(a,b') = h(a',b')$. Theorem 3 follows very quickly.
For, $H(a,b')$ is then ambient isotopic to both $H(a,b)$ and $H(a',b')$.
So $f^{a'-a}$ extends to a homeomorphism of $H(a,b)$ to itself.
Hence this power of $f$ partially extends to $H$.

Suppose the claim were not true. Then we will show by induction
that, for each positive integer $j$, there are non-negative integers
$p(j)$ and $m(j)$ and a sequence $\{ k(j,n) : n \geq s \}$
such that for each $n \geq s$,

\items
\item{(i)} $h(k(j,n), k(j,n)+p(j)) \leq 2 \ {\rm genus}(\partial H) - j$;
\item{(ii)} $0 \leq k(j,n) - n \leq m(j)$.
\enditems

\noindent This will lead to a contradiction, since $h(a,b)$ is
non-negative, but (i) implies that it is negative for
large $j$.

The induction starts with $j = 1$, where we take
$m(1) = p(1) = 0$ and $k(j,n) = n$. We now prove the
inductive step. Each interval $[k(j,n), k(j,n) + p(j)]$
corresponds to the boundary of discs. There are
a finite number of choices ($|{\cal C}|^{p(j)+1}$) for the boundary
curves of these discs. By Lemma 1,
we may assume that these discs form one of only finitely
many ($d$, say) disc arrangements.
Thus, we assign to each integer $k(j,n)$ one of these
$d$ colours. If $k(j,n')$ and $k(j,n'')$ are given the
same colour, then the pairs $(k(j,n'), k(j,n') + p(j))$
and $(k(j,n''), k(j,n'') + p(j))$ are comparable.

For any $n \geq s$, consider the values of $k(j,n')$ as
$n'$ varies between $n$ and $n + (m(j)+1) d$.
By (ii) at most $(m(j)+1)$ different $n'$ can give
the same value of $k(j,n')$. Therefore, the $k(j,n')$ take more
than $d$ different values. Therefore, we can find integers
$n'$ and $n''$ in this interval coloured with the same
colour, such that $k(j,n') \not= k(j,n'')$. Say that
$k(j,n') < k(j,n'')$.

We let 
$$\eqalign{
k(j+1,n) &= k(j,n') \cr
p(j+1) &= (m(j) +1) d + m(j) + p(j) \cr
m(j+1) &= (m(j)+1)d + m(j).}$$

We need to check that (i) and (ii) hold for $j+1$.
Since we are assuming that the claim is not true,
it must be the case that $h(k(j,n'),k(j,n'') + p(j))$ 
is strictly less than at least one of
$h(k(j,n'), k(j,n') + p(j))$ and $h(k(j,n''), k(j,n'') + p(j))$.
Hence, by (i), $h(k(j,n'),k(j,n'') + p(j)) \leq 2 \ {\rm genus}(\partial H) - j -1$. 
So,
$$\eqalign{
h(k(j+1,n), k(j+1,n)+p(j+1)) &= h(k(j,n'), k(j,n') + p(j+1)) \cr
&\leq h(k(j,n'), k(j,n'') + p(j)) \cr
&\leq 2 \ {\rm genus}(\partial H) - j -1,}$$
since
$$\eqalign{
& k(j,n') + p(j+1)  - k(j,n'') - p(j) \cr
&\qquad \geq
n' + (m(j) +1) d + m(j) + p(j)
- n'' - m(j) - p(j) \geq 0.}$$
This verifies (i) for $j+1$.

To verify (ii), note that
$$k(j+1,n) = k(j,n') \geq n' \geq n,$$
$$k(j+1,n) = k(j,n') \leq n'+ m(j) \leq n + (m(j)+1) d + m(j) = n +
m(j+1).$$
This establishes the claim, and hence completes the proof of Theorem 3.

We conclude with a proof of Theorem 4. So, let
$f \co \partial H \rightarrow \partial H$ be a
pseudo-Anosov homeomorphism. Suppose that,
for infinitely many $n$, $M \cup_{f^n \circ \phi} H$
fails to satisfy the conclusions of Theorem 1.
Then, for these $n$, $f^n\phi(C_{i(n)})$ bounds
a disc $D_n$ in $H$, for some exceptional curve $C_{i(n)}$.
We may pass to a monotone subsequence where $C_{i(n)}$
is some fixed curve $C$. But $f^n\phi(C)$ tends in $PL(\partial H)$ 
to the stable or unstable lamination of $f$, according
to whether the sequence is decreasing or increasing.
Hence, this lamination lies in the closed set
$B(\partial H)$. This completes the proof of Theorem 4.

\references

\end